\title{ALTERNATING KNOTS AND LINKS THEORY}
\author{E. PI\~NA \\
Departamento de F\'\i sica\\ Escuela Superior de F\'\i sica y Matem\'aticas del IPN \\ on sabbathical leave of absence from\\ Universidad Aut\'onoma Metropolitana - Iztapalapa, \\
P. O. Box 55 534 \ Mexico, D. F., 09340 Mexico \\
e-mail: pge@xanum.uam.mx}

\documentclass[12pt]{article}
\begin{document}
\date{ }

\maketitle

\abstract{The alternating knots, links and twists projected on the $S_2$ sphere are identified with the phase space of a Hamiltonian dynamic system of one degree of freedom. The saddles of the system correspond to the vertices, the edges correspond to the stable and unstable manifolds connecting the saddles. Each face is then oriented in one of two different senses determined by the direction of these manifolds. On a certain point of any face one could have an equilibrium elliptic point. This correspondence can be realized also between the knot (or link) and the Poincar\'e section of a two degrees of freedom integrable dynamical system. The vertices corresponding to unstable orbits (knots), and the faces to foliated torus, around a stable orbit (knot).

The associated matrix to that connected graph is decomposed in two permutations. The separation is unique for knots and is not for links. The characteristic polynomial corresponding to some knot, link or twist families was explicitly computed in terms of Chebyschev polynomials. Any vertex of a knot or link can be considered a member of one of two families, replacing the vertex by a ribbon in two different directions.}

\newpage

\section{INTRODUCTION}
Many studies on knots and links have been published in the last years with remarkable impact in the field of Physics and other sciences. One paper including many remarkable citations is the review article on the application of knot theory to Statistical Mechanics (Wu 1992). It also has been used with success in Topological Quantum Field Theory (Witten 1988, Atiyah 1990).

This study is restricted to the alternating knots that are defined when the knot is projected on a surface as a graph. One assumes over passes and underpasses alternate along the curve. The restriction is extended to links with the same property independent of the number of closed curves. In this paper we study the matrix associated to the graph and its characteristic polynomial that are important semi-invariants of the knots and links, since they can be grouped in families.

In studying knots and links several polynomials were introduced in past, like the Alexander's (1928), Jones' (1985), HOMFLY (see Freyd 1985), Kauffman's (1990), etc. (see for example Balister et al. 2001).  Some of these polynomials were written in terms of two or three variables and give therefore more information that the one is found in this paper. Nevertheless I note some simplicity and elegance in the polynomials and graphs of the approach in this paper. I expect henceforth some interest and perhaps also some future applications of this theory.

The graph of the knot or link is formed by a finite number of vertices $V$, edges $E$, and faces $F$, which obey the Euler equation
\begin{equation}
V + F = E + 2\, ,
\end{equation}
where the external face or sea is also counted as a face.
Besides, four edges join at any vertex and each edge joins two vertices, therefore one has also
\begin{equation}
E = 2 V \, , \quad F = V + 2\, .
\end{equation}

At this point of this paper I assume we have not loops, equal to faces with only one vertex or edge. Later this restriction is lifted.

If we denote by $C_j$ the number of faces having $j$ vertices and $j$ edges, then the total number of faces is
\begin{equation}
F = V + 2 = \sum_{j=2}^L C_j \, ,
\end{equation}
where $L$ is the number of vertices of the largest face of the knot.

Since edges are shared by two faces, and vertices by four faces one has
\begin{equation}
2 E = 4 V = \sum_{j=2}^L j \, C_j \, .
\end{equation}

Equating $V$ from these two equations we obtain the identity
\begin{equation}
2 C_2 + C_3 = 8 + \sum_{j=5}^L (j-4) C_j \, ,
\end{equation}
which shows the importance of the faces having two and three edges in any knot. The integer on the left is always larger than 7.

When restricted to alternating knots it is not important what lane crosses in the upper level at any vertex. There are only two possibilities for the total knot, the one we choose and the reflected one, which can be identified with the mirror image of that knot and will not be distinguished in this paper. Once one lane is identified as upper at a crossing all the next crossings are determined by the alternating restriction.

Now we introduce an orientation of all the edges of the knot that brings to mind the orientation of the stable and unstable manifolds connecting the saddles of a dynamic system in a surface. The vertex is identified with the saddle and the upper lane with the unstable manifold, leaving the two edges the vertex with opposite orientation; the down lane is associated to the stable manifold formed by two edges that are oriented pointing to the vertex. This orientation of the edges give to any face one orientation clockwise or counterclockwise, allowing the coloring of the faces of these knots or links with only two colors according to the rotation sign of the orientated boundary.

Balister et al. (Balister 2001) have considered this orientation for alternating knots or links. In their terminology it is a connected two-in two-out digraph.

To see the graphs of a simple twist and the trefoil knot in a dynamic context see for example figure 3.2.5 at page 67 of Reichl (2004).

It is clear that the assumptions here adopted about these knots or links force some relative positions between links and avoid to consider Reidemeister moves (Reidemeister 1948) of the edges of the knots and links. Besides some knots, like the square knot with $V=6$, can be presented as an alternating knot by modifying its drawing; it should be noted that they are also included in this theory.

\section{THE MATRIX OF THE ALTERNATING KNOTS}
There is an standard square matrix associated to an oriented graph (Berge 2001). The matrix determines the alternating knot and vice versa. One enumerates the vertices. The edge connecting vertices from $j$ to $k$ is represented in the matrix by the entry $M_{jk} = 1$ at row $j$, column $k$ of the associated matrix. Exceptionally in some links two edges connect the same two vertices. In that cases the entry corresponding to those vertices will be equal to 2; the other entries are zero. The resulting matrix is an invariant of the knot, it can be think as equivalent to the knot. For Balister et al. (Balister 2001) it is the {\sl adjacency} matrix. In what follows one will present some properties of this matrix.

From this alternating matrix we can identify if it is a knot or a link and how many different closed paths form it. We form subgraphs starting at any non zero entry $M_{jk}$ of the associated matrix, add the edge in the same row represented by the non zero entry $M_{jl}$, then we add the edge in the same column $M_{ml}$, then the edge in the same row $M_{mn}$ and so successively until the subgraph is closed with an edge $M_{zk}$ in the same column $k$ that the starting element. The resulting subgraph correspond to a knot or an unknotted element. A link formed by $N$ elements has $N$ subgraphs constructed in this way. Note that an entry 2 of the matrix correspond to a circle formed by two edges and appears only in links; by itself represent one component of the link.

The algebraic version of the above properties is that the matrix of a knot is the sum of two permutation matrices, that are identified in alternated positions in the subgraph. This decomposition is unique for a knot. Each element of a link has a unique decomposition, but now the matrix is not decomposed in a unique form since the decompositions corresponding to each component in the link are independent.

Another property is that any column and any row of the matrix has two 1's (or only one 2) and therefore the vector that has all the entries with the same value is an eigenvector of this matrix with eigenvalue 2. All the other eigenvectors are orthogonal to this vector, that means that the sum of all its components is null.

The characteristic polynomial of this matrix is a polynomial with integer coefficients since all the components of the matrix are integers. The roots of the polynomial are algebraic numbers. The polynomial factors always in two or more polynomials having also integer coefficients. Factor $(x - 2)$ appears in all the polynomials.

The characteristic polynomial of the matrix of a knot or link is of degree $V$. The coefficient of the power $V-1$ is zero since we have not loops, nor faces of one vertex and of one edge. The coefficient of power $V-2$ is $-C_2$. The coefficient of power $V-3$ is $-C_3$. Not simple rule exist for other coefficients. The proof of these properties comes from the relation between the coefficients of the characteristic polynomial and the traces of the powers of the matrix. One also uses that the trace of the original matrix is always zero. All this comes from the well known property that the number of closed paths of length $k$ in a graph are equal to the trace of the power $k$ of the matrix (Andrews 1999).

How good is this polynomial to identify a particular knot or link? The answer to this important question is not always positive. There are different links with same polynomial. And the same link with different positions of the components produces sometimes different polynomials. Nevertheless this polynomial is an important invariant of the alternating knots and links.

Many knots have a cyclic symmetry that is also present in the matrix and their characteristic polynomial as a representation of the symmetry. For example the torus knots of $2 N + 1$ vertices have an associated matrix that can be represented as the sum of a cyclic permutation plus the inverse permutation. The power $2N+1$ of these permutations is the identity matrix. It follows that this matrix has a minimal polynomial related to the group. However in this study this family is one half of a larger family that will be presented later in this paper.

There are other symmetries that imply some factorization of the characteristic polynomials of a knot or link. For example some of these matrices commute with a permutation self inverse matrix; this property produces that the subspaces spanned by the eigenvectors associated to projectors related to that permutation matrix can be factorized in a number of polynomials having as roots the corresponding eigenvalues. Or some matrices of this representation of knots and links can be separated in a finite number of sub matrices of lower dimension. Also implying a factorization of the characteristic polynomial.

\section{FAMILIES OF LINKS AND KNOTS}
The matrices of many knots and links were studied and the characteristic polynomials of these matrices were computed and compared.

Many knots have an isolated face of two edges or a series of contiguous faces of two edges, forming a {\sl ribbon}. One notes that one face of two edges can be destroyed by suppressing the two edges and colliding or joining the two vertices in one vertex without modifying the rest of the knot or link and the orientation of the edges forming it. This operation can be iterated in a series of contiguous faces of two edges (ribbon), but also can be reversed to starting with a knot and obtaining other knot with an extra vertex in the same direction of the ribbon. This operation gives birth to a family of knots having the same structure but different number of faces of two edges in a contiguous series of a ribbon. The knots or links form a family and the characteristic polynomial of any member of the family can be computed as a function of three contiguous polynomials. Looking to the simplest families of polynomials one finds useful to follow the destruction of a series and suppress also the last vertex by joining two times one incoming edge to an outgoing edge, the gap in the same direction of the ribbon. This is the so called (Farmer 1995) {\sl crossing elimination}, upon choosing from the two possibilities the ribbon direction. This produces in most cases other knot or link, but also a twist characterized by having an unknotted component with one or two loops, that could be unknotted. Also when we can separate two parts of a knot (link) by two edges, a {\sl composed} knot (Farmer 1995), we can twist these two edges and form a family of twist than could be unknotted by the characteristic Reidemeister move.

With this construction any vertex can bifurcate in two families, with two possible orientations for the family. Or any vertex of a knot or link without faces of two loops can be destroyed by thinking it as the last member of the family. These operations remind the Skein relation for the Kauffman polynomial (Kauffman 1990), other kind of polynomials associated to knots and links.

As an example consider the Fullerene graph (Hartsfield 2003) of 60 carbon atoms that has double carbon bonds on three alternated edges of contiguous hexagons, and a simple bond separating each pair hexagon-pentagon. This can be regarded as an alternating graph of 60 vertices with icosahedral symmetry. Applying the previous operation the pair of carbon atoms with a double bond is replaced by only one carbon atom, without double bonds; it comes a molecule with 30 carbon atoms, preserving the symmetry.

The matrix of a twist can be computed in the same way we do for any knot or link. The only difference is that in some cases we have now loops: one or two faces in a link with the property of having only one edge. The polynomials of a family of twists are simpler that other family polynomials and we start with them. These are used afterward, in some cases, as the first polynomial of a family.

We consider the simplest family of twist formed by twisting a circle. The first twist produces a figure eight with two loops. The next twist in this family has two loops joined to a face of two edges. The general term of the family is a twist with $V$ vertices, two loops, $V-1$ faces of two edges and the external face of $2V$ edges. The polynomials of this family obeys the recurrence relation that is found by induction
\begin{equation}
P_{V+2}(x) - x P_{V+1}(x) + P_V(x) = 0\, ,
\end{equation}
and can be determined by the two first. All these polynomials have the factor $x-2$ and are expressed as
\begin{equation}
P_V(x) = (x - 2) J_{V-1}(x)
\end{equation}
in terms of Polynomials $J_k(x)$ which are well known polynomials, particular cases of Jacobi, Gegenbauer and Chebyschev polynomials denoted as
\begin{equation}
J_{k}(x) = U_{k}(x/2)\, ,
\end{equation}
obeying the recurrence relation (6).
They can be written explicitly (Bell 2004) as
\begin{equation}
J_k(x) = \sum_{j=0}^{[k/2]} (-1)^j \frac{(k-j)!}{j! (k-2j)!} x^{k-2j}\, ,
\end{equation}
with the generating function
\begin{equation}
\frac{1}{1 - t x + t^2} = \sum_{k=0}^\infty t^k \, J_{k}(x)\, .
\end{equation}

The first polynomials of this family are $1$, $x$, $x^2-1$, $x^3-2x$, $x^4-3x^2+1$, $x^5-4x^3+3x$, $x^6-5x^4+6x^2-1$, $x^7-6x^5+10x^3-4x$. These polynomials will be very important in this study of alternating knots and links.

A useful quadratic relation between two contiguous polynomials is
$$
J_k^2(x) - x J_k(x) J_{k-1}(x) + J_{k-1}^2(x) = 1\, .
$$

A second family of twists is formed starting from the Hopf link of two circles and twisting one of the 4 edges of this link.
The polynomials corresponding to this family of twists obeys the same recurrence relation (6) with different initial conditions. The polynomial of $V$ vertices of this graph is
\begin{equation}
(x-2)[(x + 2) J_{V-2}(x) - x J_{V-3}(x)]\, ,
\end{equation}
in terms also of the Chebyshev polynomials.

The polynomials of the twist family of the trefoil knot obey the same recurrence relation (6) with other initial conditions that give the family of polynomials
\begin{equation}
(x-2)(x+1)[(x + 1) J_{V - 3}(x) - x J_{V - 4}(x)]\, .
\end{equation}

The polynomials of the twist family formed by twisting one edge of the knot of 4 vertices is similarly constructed starting with the two polynomials $x^4 -2 x^2 - 4 x$ and $x^5 - x^4 - 2 x^3 + x - 2$. that determine the family.

We have also the twist families formed by twisting two edges that separate in two parts a knot or link. The knot with these two parts is the product of knots named composition (Farmer 1995) that give birth to the concept of {\sl prime knot} settled by H. Schubert in 1949. The polynomials of each family of twist satisfies the recurrence (6) that is a characteristic of the twists.

These twist families were used by connecting one member of one of these families to the part where it originates, by colliding in one new vertex that produces the loop to be converted in a two edges face, the twist becoming in this way a knot or a link. The necessary condition for allowing a collision is that the two edges to be connected, the one of the loop with another should have an antiparallel orientation with respect to the loop orientation. The approach of the tail having at the extreme a loop to the other part of the twist can be thought as a motion over the surface where the knot is drawn and this transformation does not modify the condition of being an alternating knot as was previously defined. But a different possibility presents itself if we think this tail as an scorpion tail that can connect the loop with an internal edge of the twist. This results in a knot or link that is not trivially an alternating knot or link, since perhaps cannot be drawn on the same surface where the twist was constructed. However it is an alternating knot (link) by construction and this theory is still useful for this knot (link).

This last remark has not been developed in this paper because it modifies the implicit hypothesis of having the graph of the knot (link) on a sphere or plane. By modifying the genus of the surface one can consider other alternating knots but also the knot could lose their crossings. Then we abandon this possibility in this paper.

The polynomials of knots (links) are going to form a family by taking as the two first polynomials the two corresponding ones to the collision of two antiparallel oriented edges before and after collision, plus the next members of the family that are formed by creating a new knot (link) by adding a two edge face in parallel direction to the initial antiparallel direction of collision.

The polynomials of the family will obey now a recurrence relation, that is different to the twist recurrence (6) by the fact that it has now a non-homogeneous term. But since all the polynomials have the factor $x-2$ it should be present the same factor in the source and the recurrence is always of the form
\begin{equation}
P_{V+2}(x) - x P_{V+1}(x) + P_V(x) = (x-2) H(x)\, ,
\end{equation}
where H(x) is a typical polynomial of the family. This recurrence can be rewritten by transforming to the homogeneous form (6) by adding to any polynomial of the family the factor $H(x)$ of the source as
\begin{equation}
P_{V+2}(x) + H(x) - x [P_{V+1}(x) + H(x)] + [P_V(x) + H(x)] = 0 \, ,
\end{equation}
that can be solved as the twist examples. The polynomial $P_V(x)$ of a knot (link) is a linear combination of the Chebyshev polynomials minus the $H(x)$ polynomial.

The classification used for knots follows the order presented by Rolfson 1976) for prime knots and links, the first number denotes the number of vertices, the subindex is the place in the table of prime knots published in the Appendix C of the Rolfsen's book. This classification is the Standard one. For links the super index denotes the number of components of the link in Appendix C of the same book.

As an example we consider the family of cyclic torus knots (links) formed by a chain of $V$ faces of two edges and two other faces of $V$ vertices (edges). When $V$ is an odd number it is a knot, when $V$ is an even number it is a link formed by two interwoven circles. The family begins with the eight shape twist of $V=1$, then the Hopf link, the trefoil knot, the Solomon's knot, $5_1$, $6_1^2$, $7_1$, $8_1^2$, etc. The polynomials of the knots (links) of this family are of the form
\begin{equation}
2 [J_{V}(x) - 1] - x J_{V-1}(x)\, ,
\end{equation}
showing $H(x)=2$ for this family. But this expression (15) hides many other symmetry features of this family of polynomials. For $V$ an odd number we have a knot with polynomial
\begin{equation}
2 [J_{2 k + 1}(x)-1] - x J_{2 k} = (x - 2) [J_k(x) + J_{k - 1}(x)]^2\, .
\end{equation}
For $V$ an even number we have a link with
\begin{equation}
2 [J_{2 k}(x)-1] - x J_{2 k-1} = (x^2 - 4) [J_{k-1}(x)]^2\, .
\end{equation}
The square factor in these expressions repeats itself as a single factor in the polynomial of many other knots and links.

The polynomials of the family including the twist of two vertices, the trefoil, the knot of 4 vertices, the knot $5_2$, the knot $6_1$, the knot $7_2$, the knot $8_1$, the knot $9_2$, the knot $10_1$, etc., called twist knots (Farmer 1995), are
\begin{equation}
(x^3 -x - 2) J_{V-3}(x) - x^2 J_{V-4}(x) - 2 x\, ,
\end{equation}
where $V$ are the number of vertices of any member of this family. Polynomial $H(x)$ is $2 x$ in this case. 

Similar families are the one beginning with the simple twist of 2 loops and 2 faces of two edges, followed by the family of links and knots $4_1^2$, $5_2$, $6_2^2$, $7_3$, $8_2^2$, $9_3$ etc.

This three families begin with a link and the family is produced by collision of two loops of this 
link. Each family has a ribbon. The three can be considered as a family of two ribbons.
To families formed by two ribbons belongs a knot with $k - 1$ faces of two edges connected to $j - 1$ faces of two edges, limited by $k$ vertices and $j$ vertices. It is a knot of $V=k+j$ vertices, having also two faces of $k+1$ vertices and two faces of $j+1$ vertices. It has a polynomial $F_{j,k}(x)$ given by
\begin{equation}
F_{j, k}(x) = J_k(x) J_j(x) - J_{k-2}(x) J_{j-2} - 2 J_{j-1}(x) - 2 J_{k-1}(x)\, ,
\end{equation}
that is symmetric with respect to interchange of the two indexes $j$ and $k$. Examples of this family are the cyclic polynomials (15), when $k=1$, $V=j+1$, and also the polynomials corresponding to the knots $F_{2,2} \Leftrightarrow 4_1$, $F_{3,2} \Leftrightarrow 5_2$, $F_{4,2} \Leftrightarrow 6_1$, $F_{5,2} \Leftrightarrow 7_2$, $F_{4,3} \Leftrightarrow 7_3$, $F_{6,2} \Leftrightarrow 8_1$, $F_{4,4} \Leftrightarrow 8_3$, $F_{7,2} \Leftrightarrow 9_2$, $F_{6,3} \Leftrightarrow 9_3$, $F_{5,4} \Leftrightarrow 9_4$, $F_{8,2} \Leftrightarrow 10_1$, etc. And the links $F_{3,3} \Leftrightarrow 6_2^2$, $F_{5,3} \Leftrightarrow 8_2^2$, etc.

An important example is obtained if we consider two different links. One is formed by three circles connected as links of a simple chain with $V=4$. The other is the cyclic four vertices link of two circles $4_1^2$ called Solomon's knot (Farmer 1995). They can collide each to form the same link of 5 vertices $5_1^2$ formed by an eight shaped twist connected to a circle. The two families are different, growing in orthogonal directions and having some different components. But they have the same family of polynomials
\begin{equation}
x^2 [(x^2 - 2) J_{V-4}(x) - 2 J_{V-5}(x) - 2]\, ,
\end{equation}
where we need $J_{-1}(x)=0$. The first family is $5_1^2$, $6_3^2$, $7_3^2$, $8_6^2$, $9_{10}^2$, etc. The second family is $5_1^2$, $6_3^3$, $7_3^2$, etc.

A new family starts with the first twist of the trefoil knot followed by $5_2$, $6_2$, $7_4$, $8_4$, $9_6$, $10_4$, etc.

Other families begin with one knot and the family is produced by collision of two edges of that knot. One example begins with the $5_1$ knot followed by $6_2$, $7_5$, $8_6$, $9_7$, $10_{20}$, etc. We note that any knot belongs to different families.

These two families are particular cases of a more general family.
Unstinted use of induction lead us to the three ribbon knot or link with $V=k+l+m$, $k+l+m-3$ faces of two edges, two faces of $m+2$ edges, one face of $k+l$ edges, one face of $k+1$ edges, and one face of $l+1$ edges. The corresponding polynomial is
$$
P_{k,l;m}(x) = [J_{k-2}(x) J_{l-2}(x) + J_k(x) J_l(x) - 2] J_m(x)
$$
\begin{equation}
- x [J_{k-1}(x) +  J_{l-1}(x) +  J_{k-2}(x) J_{l-2}(x) - 1] J_{m-1}(x) - 2 J_{k-1}(x) J_{l-1}(x)\, ,
\end{equation}
which is symmetric in the $k$ and $l$ indices. Some Examples are $P_{3,2;1} \Leftrightarrow 6_2$, $P_{3,3;1} \Leftrightarrow 7_4$, $P_{3,2;2} \Leftrightarrow 7_5$, $P_{5,2;1} \Leftrightarrow 8_2$, $P_{4,3;1} \Leftrightarrow 8_4$, $P_{3,2;3} \Leftrightarrow 8.6$, $P_{5,2;2} \Leftrightarrow 9_6$, $P_{3,2;4} \Leftrightarrow 9_7$, $P_{4,3;2} \Leftrightarrow 9_9$, $P_{3,3;3} \Leftrightarrow 9_{10}$, $P_{6,1;3} \Leftrightarrow 10_4$, $P_{5,2;3} \Leftrightarrow 10_6$, $P_{4,3;3} \Leftrightarrow 10_{11}$, $P_{5,2;3} \Leftrightarrow 10_{20}$, $P_{2,2;2} \Leftrightarrow 6_2^3$, $P_{2,2;3} \Leftrightarrow 7_2^3$, $P_{4,2;2} \Leftrightarrow 8_2^3$, $P_{3,3;2} \Leftrightarrow 8_2^4$.

A similar knot formed by three ribbons with $V=j+k+l$ vertices, one face of $k+l$ edges, one face of $l+m$ edges, one of $m+k$ edges, and two faces of three edges leads to the symmetric polynomial in the three indexes
$$
G_{k,l,m}(x) = x[J_{k-1}(x)J_l(x)J_m(x)+J_k(x)J_{l-1}(x)J_m(x)+J_k(x)J_l(x)J_{m-1}(x)]
$$
$$
-x^2[J_k(x)J_{l-1}(x) J_{m-1}(x)+J_{k-1}(x)J_l(x) J_{m-1}(x)+J_{k-1}(x) J_{l-1}(x)J_m(x)]
$$
\begin{equation}
+ (x^3-2)J_{k-1}(x)J_{l-1}(x)J_{m-1}(x)-x[J_{k-1}(x)+J_{l-1}(x)+J_{m-1}(x)]\, .
\end{equation}
We have the property $P_{k,l;1}=G_{k,l,1}$.  Many links and knots share this family. Some typical examples are $G_{1,1,1} \Leftrightarrow 3_1$, $G_{3,3,2} \Leftrightarrow 8_5$, $G_{3,3,3} \Leftrightarrow 9_{35}$, $G_{5,3,2} \Leftrightarrow 10_{46}$, $G_{4,3,3} \Leftrightarrow 10_{61}$, $P_{2,2;1} \Leftrightarrow 5_2^1$, $G_{4,2,1} \Leftrightarrow 7_1^2$, $G_{3,2,2} \Leftrightarrow 7_4^2$, $G_{4,2,2} \Leftrightarrow 8_1^3$. The composition of two cyclic polynomials of $k$ and $l$ vertices have the polynomial $G_{k,l,0}$, where we assume again $J_{-1}(x) = 0$.

The particular case when all the indices are the same $k=l=m$, $(V=3 k)$, has the cyclic symmetry of the trefoil knot. The polynomial simplifies to the form
\begin{equation}
G_{k,k,k}(x) = (x - 2)(1 + x)^2 J_{k-1}(x)^3 \, , 
\end{equation}
where the factor $(x - 2)(1 + x)^2$ corresponds to the polynomial of the trefoil knot.

A closed simple chain of $V=2 k$ vertices produces the polynomial
\begin{equation}
P_k(x) x^k\ ,
\end{equation}
where $P_k(x)$ is the cyclic polynomial (15) of the cyclic torus knots (links).

The generalization of these two families is formed by the family of $k$ ribbons, all of the same length of $m$ vertices, $V= k m$, with a cyclic $k$-symmetry and with polynomial
\begin{equation}
P_k(x) J_{m-1}^k(x)\, .
\end{equation}

Last in our examples we found the chaining of one cyclic knot of $n$ vertices with a simple chain of $k$ links $(V = 2 n + 2 k)$:
$$
L_{k;n}(x) = x^k [2 J_k(x)- x J_{k-1}(x)][J_n(x) - 1] + 
$$
\begin{equation}
x^k [- x J_k + (x^2 - 2)J_{k-1}(x)] J_{n-1}(x)\, .
\end{equation}

These are the simplest cases. A family with less symmetry begins with $5_2$, and then $6_3$, $7_6$, $8_8$, etc. Therefore any alternating knot (link) can be considered a member of a family where any vertex not belonging to a ribbon can be replaced by a ribbon. A classification of knots (links) follows according to the number of ribbons and its relative connections, each ribbon connected by four different edges to the rest of the knot (link). Avoiding composition (and therefore twists also), it is possible to start a classification schema. The family of one ribbon is the torus cyclic family. The two ribbons family include the twist knots (links). The three ribbons families are two and for both were computed the corresponding polynomial. The four ribbons families are five, from these only some particular cases have been computed. 

I present this part of the theory assuming it is interesting without answering many questions that could be addressed in connection with this perspective. 

\

{\large \bf ACKNOWLEDGMENTS }

\

I thank to Mathematician E. Virue\~na for teaching me the algorithm to draw geometric pentagonal patterns of Moorish art, that was at the origin of this study; and to Profs. Mar\'\i a Teresa de la Selva and Lidia Jim\'enez-Lara for helping comments.

\

{\large \bf REFERENCES }

\

\noindent G. E. Andrews, R. Askey, and R. Roy {\sl Special Functions}, (Cambridge University Press, Cambridge, 1999).

\

\noindent J. W. Alexander "Topological Invariants of Knots and Links" {\sl Trans. Amer. Math. Soc.} {\bf 30} 275-306 (1928).

\

\noindent M. F. Atiyah {\sl The Geometry and Physics of Knots}, (Cambridge University Press, Cambridge, 1990).

\

\noindent P. N. Balister, B. Bollob\'as, O. M. Riordan and A. D. Scott "Alternating Knot Diagrams, Euler Circuits and the Interlace Polynomial" {\sl Europ. J. Combinatorics} {\bf 22} 1-4 (2001).

\

\noindent W. W. Bell {\sl Special Functions for Scientist and Engineers}, (Dover, New York, 2004).

\

\noindent C. Berge {\sl The Theory of Graphs} (Dover, New York, 2001).

\

\noindent David. W. Farmer and T. B. Stanford {\sl Knots and Surfaces. A guide to Discovering} (American Mathematical Society, Providence, 1995).

\

\noindent P. Freyd, D. Yetter, J. Hoste, W. B. R. Likorish, K. Millet, and A. Ocneanu (HOMFLY) "A New Polynomial Invariant of Knots and Links" {\sl Bull. of the Am. Math. Soc.} {\bf 12} 239-246
(1985) .

\

\noindent N. Hartsfield, G. Ringel {\sl Pearls in Graph Theory} (Dover, Mineola, 2003) p. 22.

\

\noindent V. F. R. Jones "A polynomial invariant for links via von Neumann algebra"{\sl Bull. of the Am. Math. Soc.} {\bf 12} 103-111 (1985).

\

\noindent L. H. Kauffman "An invariant of regular isotopy" {\sl Trans. Amer. Math. Soc.}, {\bf 318} 417-471 (1990).

\

\noindent L. E. Reichl {\sl The Transition to Chaos} 2nd ed. (Springer, 
New York, 2004). p 67.

\

\noindent K. Reidemeister {\sl Knotentheorie} (Chelsea, New York, 1948).

\

\noindent D. Rolfsen {\sl Knots and Links}, (Publish or Perish, Berkeley, 1976).

\

\noindent E. Witten "Topological quantum field theory" {\sl Commun. Math. Phys.} {\bf 117} 353-386 (1988).

\

\noindent F. Y. Wu "Knot theory and statistical mechanics"
{\sl Reviews of Modern Physics} {\bf 64} 1099-1131 (1992). 

\end{document}